\documentclass[a4paper]{article}

\usepackage{authblk}
\usepackage{amsfonts}
\usepackage{amsmath}
\usepackage{amsthm}
\usepackage{amssymb}
\usepackage{mathrsfs}
\usepackage{color}
\usepackage[table]{xcolor}

\usepackage[pagewise]{lineno}

\numberwithin{equation}{section}

\newtheorem{Theorem}{Theorem}  
\newtheorem{Lemma}{Lemma}

\newenvironment{Proof}
{\noindent
\textsc{Proof -- }}{\hfill$\square$\medskip}

\newcounter{remark} 
\renewcommand{\theremark}{\arabic{remark}}
\newenvironment{Remark}{
\refstepcounter{remark}
\bigskip\noindent \textit{Remark
\theremark\ - }}

\begin{document}

\title{$L^2$-setting theory for the solutions to\\ 2D Navier--Stokes equations: some new estimates
}
\author{Alfonsina Tartaglione\footnote{Dipartimento di Matematica e Fisica, Universit\`a degli Studi della Campania ``Luigi Vanvitelli", viale Lincoln 5, 81100 Caserta (Italy); e--mail: alfonsina.tartaglione@unicampania.it\\
\\
{\bf Conflict of Interest Statement.} The author declares there is no conflict of interest.\\
{\bf Data Availability Statement.} Data sharing not applicable to this article as no datasets were generated or analysed during the current study.
}
}
\date{}
\maketitle

\begin{abstract}
For an unbounded planar domain $\Omega$, we prove that the solution $u$ to the IBVP for the Navier--Stokes equations with initial datum $u_0\in L^2(\Omega)$ satisfies the following estimate 
    \begin{equation*}
        \int_0^{+\infty}\|u(\tau)\|_\infty^2d\tau\le c(1+\|u_0\|_2^2)\|u_0\|_2^2,
    \end{equation*}
    established by R. Farwig and Y. Giga [{\sl Algebra i Analiz}, 36, 289-307 (2024)] for bounded domains. Moreover, we show that
    \begin{equation*}
    \int_0^{+\infty}\tau\|P\Delta u(\tau)\|_2^2d\tau\le c(1+\|u_0\|_2^4 e^{c_1\|u_0\|_2^4})\|u_0\|_2^2,
\end{equation*}
and
\begin{equation*}
    \int_0^{+\infty}\tau\|u(\tau)\|_\infty^4d\tau\le c(1+\|u_0\|_2^4 e^{c_1\|u_0\|_2^4})\|u_0\|_2^4.
\end{equation*}
\end{abstract}

\section{Introduction}

Let us consider the motions of an incompressible fluid 
in a domain $\Omega$ of ${\mathbb R}^2$ with $C^2$ boundary. In this paper,

\begin{equation}\label{dominio}
\Omega\equiv
\begin{cases}
\text{a bounded } C^2\text{-domain};\\
\text{an exterior } C^2\text{-domain};\\
\mathbb{R}^2;\\
\mathbb{R}^2_+.
\end{cases}
\end{equation}

The pair $(u,p)$, where $u$ and $p$ denote the velocity and pressure fields, respectively, satisfies the Navier--Stokes equations
\begin{equation}\label{NavierStokes}
\begin{array}{r@{}l}
    u_t-\Delta u + u \cdot \nabla u + \nabla p &{} = 0\qquad {\mathrm {in }}\ (0,+\infty)\times\Omega,\\
    \nabla \cdot u &{} = 0\qquad {\mathrm {in }}\ (0,+\infty)\times\Omega.
\end{array}
\end{equation}
We complement \eqref{NavierStokes} with the initial and boundary conditions
\begin{align}\label{condition1}
    u_{|t=0}=u_0 \quad &{\hbox {in}}\ \Omega,\\\label{condition2}
    u=0 \quad&{\hbox{on}}\ (0,+\infty)\times\partial\Omega.
\end{align}
Denoting by $J^2(\Omega)$ the completion of the spatial vector fields belonging to $C_0^\infty(\Omega)$ and with free divergence, it is well--known that, for all $u_0\in J^2(\Omega)$, there exists a unique (weak) solution $(u,p)$ 
to \eqref{NavierStokes}--\eqref{condition2}, such that, for all $t>0$, $u$ satisfies the energy identity
\begin{align}\label{identity}
    \frac{1}{2}\int_\Omega|u(t,x)|^2dx+\int_0^t\int_\Omega|\nabla u(\tau,x)|^2dxd\tau=\frac{1}{2}\int_\Omega |u_0(x)|^2dx.
\end{align}

In a recent paper, R. Farwig and Y. Giga \cite{Farwig} proved that, in the case of bounded domains, the solution $u$ also satisfies the following estimate:
\begin{equation}\label{estimate}
        \int_0^{+\infty}\|u(\tau)\|_\infty^2d\tau\le c(1+\|u_0\|_2^2)\|u_0\|_2^2,
    \end{equation}
 which cannot be deduced directly  from \eqref{identity} via the Sobolev embedding, since the involved 
exponent is critical ($p=2$). One of the purposes of this paper is to extend estimate \eqref{estimate} to the case of the other domains satisfying \eqref{dominio}. Precisely, we are going to prove the following

\begin{Theorem}\label{maintheorem}
Let $\Omega$ be as in \eqref{dominio} and $(u,p)$ be the solution to \eqref{NavierStokes}--\eqref{condition2}, with $u_0\in J^2(\Omega)$. Then, $u$ satisfies 
\begin{equation}\label{estimateripetuta}
        \int_0^{+\infty}\|u(\tau)\|_\infty^2{\hbox {d}}\tau\le c(1+\|u_0\|_2^2)\|u_0\|_2^2,
    \end{equation}
    where $c$ is a constant independent of $u_0$.
\end{Theorem}

To introduce the result stated in Theorem \ref{maintheorem}, let us consider the Stokes system obtained by neglecting the nonlinear term in \eqref{NavierStokes}: 
\begin{equation}\label{Stokes}
\begin{array}{r@{}l}
    w_t-\Delta w + \nabla q &{}=0\qquad {\mathrm {in }}\ (0,+\infty)\times\Omega,\\
    \nabla \cdot w&{}= 0\qquad {\mathrm {in }}\ (0,+\infty)\times\Omega,
\end{array}
\end{equation}
with the same initial and boundary conditions: 
\begin{align}\label{condition7}
    w_{|t=0}=u_0 \quad&{\hbox {in}}\ \Omega,\\\label{condition8}
    w=0 \quad&{\hbox{on}}\ (0,+\infty)\times\partial\Omega.
\end{align}
Observe that the solution $(w,q)$ to \eqref{Stokes}--\eqref{condition8} satisfies the same energy identity \eqref{identity}. This is a consequence of the incompressibility and boundary conditions, which ensure that the nonlinear term in \eqref{NavierStokes} does not contribute to the energy balance.

Set $v=u-w$ and $\pi=p-q$. Of course, the pair $(v,\pi)$ solves the problem

\begin{equation}\label{forceterm}
\begin{array}{r@{}l}
    v_t-\Delta v + \nabla \pi  =-u\cdot\nabla u\qquad\qquad &{\mathrm {in }}\ (0,+\infty)\times\Omega,\\
    \nabla \cdot v = 0\qquad\qquad\qquad &{\mathrm {in }}\ (0,+\infty)\times\Omega,\\
\end{array}
\end{equation}
with zero initial and boundary conditions: 
\begin{align}\label{nullcondition1}
    v_{|t=0}=0 \quad{\hbox {in}}\ &\Omega,\\\label{nullcondition2}
    v=0 \quad{\hbox{on}}\ &(0,+\infty)\times\partial\Omega.
\end{align}

In \cite{Farwig}, estimate \eqref{estimate} for bounded domains is established  by proving
\begin{equation}\label{stimalineare}
    \displaystyle\int_0^{+\infty}\|w(\tau)\|_\infty^2{\text{d}}\tau \le c\|u_0\|_2^2,
\end{equation}
and
\begin{equation}\label{stimanonlineare}
    \displaystyle\int_0^{+\infty}\|v(\tau)\|_\infty^2{\text{d}}\tau \le c\|u_0\|_2^4.
\end{equation}

We first note that \eqref{stimalineare} in \cite{Farwig} is derived using a generalized Marcinkiewicz interpolation theorem based on elementary Stokes semigroup theory. This argument carries over directly to the cases of exterior domains, the half--plane and the whole plane. As for \eqref{stimanonlineare}, we provide a new proof that relies on different techniques from those employed in \cite{Farwig}. 
Specifically, to prove Theorem \ref{maintheorem} we exploit properties of the solutions to the Stokes problem with initial data in $L^1(\Omega)$, as developed in \cite{Mare2, Mare3} and their application in \cite{Mare4}, and apply an 
$L^1-L^\infty$ duality argument. 
Of course, we do not exclude that the tools of \cite{Farwig} lead to our same result, but we think that our proof is not only very short, but also does not require any additional background theory.

\medskip

We go beyond the result stated in Theorem \ref{maintheorem}, by investigating the time integrability of the square norm of the second derivatives of the solutions to the Navier--Stokes equations. We establish a new estimate that seems to be currently unknown already for bounded domains, and for the solutions to the Stokes system, too.
By the semigroup properties of the Stokes operator, it is well--known that 
\begin{equation}\label{comportamentoneltempo}
\|P\Delta u(t)\|_2\le c\|u_0\|_2 t^{-1}.
\end{equation}
The behavior \eqref{comportamentoneltempo} does not allow one to state that $\|P\Delta u\|_2\in L^2(0,T),\ \forall T>0$, but only that $\|P\Delta u\|_2\in L^2_{\text{weak}}(0,T),\ \forall T>0$. The same conclusion holds for $\|u\|_\infty^2$, taking into account that, in virtue of the following inequality (see \cite{Mare3})
\begin{equation}\label{FKinequality1}
    \|u(t)\|_\infty\le c\|P\Delta u(t)\|_2^{1/2}\|u_0\|_2^{1/2}, \quad\forall t>0,
\end{equation}
one has
\begin{equation}\label{comportamentoneltempo2}
\|u(t)\|_\infty^2\le c\|u_0\|_2^2 t^{-1}.
\end{equation}

So, we consider as interesting the following result.

\begin{Theorem}\label{derivateseconde}
Let $\Omega$ be as in \eqref{dominio} and  $(u,p)$ be the solution to \eqref{NavierStokes}--\eqref{condition2}, with $u_0\in J^2(\Omega)$. Then, $u$ satisfies
\begin{equation}\label{primastimateorema2}
    \int_0^{+\infty}\tau\|P\Delta u(\tau)\|_2^2d\tau\le c(1+\|u_0\|_2^4 e^{c_1\|u_0\|_2^4})\|u_0\|_2^2,
\end{equation}
and
\begin{equation}\label{secondastimateorema2}
    \int_0^{+\infty}\tau\|u(\tau)\|_\infty^4d\tau\le c(1+\|u_0\|_2^4 e^{c_1\|u_0\|_2^4})\|u_0\|_2^4,
\end{equation}
where $c$ and $c_1$ are constants independent of $u_0$.
\end{Theorem}
This is a surprising result, as we introduced a weight that worsens the behavior at infinity of $\|P\Delta u\|_2^2$ or $\|u\|_\infty^4$ while not improving their behavior near zero. Actually, this appears in accordance with estimate \eqref{estimateripetuta}. 

\medskip

    We emphasize that estimate \eqref{primastimateorema2} seems to be new even for bounded domains, as well as for  the solutions to the Stokes system, while in the linear case, or when the domain is bounded, \eqref{secondastimateorema2} can be deduced from \eqref{estimateripetuta} and \eqref{comportamentoneltempo2}. 
We show that, for more general domains, neither \eqref{estimateripetuta} or \eqref{comportamentoneltempo2} is required. The main purpose of this note is just to highlight a sort of congeniality of the $L^2$-theory for {\underline{2D Navier--Stokes equations}, meaning that some inequalities available in the linear theory can be fully transferred to the nonlinear setting, provided that the right-hand side is expressed through a suitable combination of the $L^2$-norm of the initial datum.

\begin{Remark}
We point out that the analysis performed in this paper is restricted to the two--dimensional case, as this allows us to avoid imposing an additional Prodi--Serrin type condition \cite{Prodi, Serrin}, which would be 
required in dimension $n\ge 3$ (as for estimate \eqref{estimateripetuta} cfr. \cite{Farwig} in the case of bounded domains).
\end{Remark}

\begin{Remark}
We observe that our approach, with suitable modifications and concerning the possible singularities in $0$, also applies to the solutions to the nonlinear equation 
\begin{equation}
u_t-\nabla\cdot(|\nabla u|^{p-2}\nabla u)=0,\quad p\in(1,2),
\end{equation}
in the wake of the results in \cite{CrispoDiFeola} and \cite{DiFeolaRuzicka}. 
\end{Remark}

\bigskip

The paper is organized as follows. In Section \ref{preliminary}, we collect some auxiliary results. Sections \ref{sectionproof} and  \ref{sectionproof2} contain the proofs of Theorem \ref{maintheorem} and Theorem \ref{derivateseconde}, respectively.\footnote{The present note is an improvement of the paper \cite{Tartaglione}.}

\section{Notation and auxiliary results}\label{preliminary}

Throughout the paper we denote by $\Omega$ a domain of ${\mathbb R}^2$ with $C^2$ boundary. $\Omega$ may be of the type described in (1). 
We use the standard notation $L^p(\Omega)$ and $W^{m,p}(\Omega)$  
for the Lebesgue and Sobolev spaces, respectively, and $\|\cdot\|_p$, $\|\cdot\|_{m,p}$ for the corresponding norms. 
Moreover,

\medskip

 $\bullet\ {\mathscr{C}}_0(\Omega):=\{\psi: \psi\in C_0^\infty(\Omega)\ {\mathrm{and}}\ \nabla\cdot\psi=0\}$,

\medskip 

 $\bullet\ J^\infty(\Omega):= {\text{completion of }}{\mathscr{C}}_0(\Omega) {\text{ in the }} L^\infty(\Omega){\text{ norm}}$,

\smallskip 

 $\bullet\ J^p(\Omega):={\text{completion of }}{\mathscr{C}}_0(\Omega) {\text{ in the }} L^p(\Omega){\text{ norm}}$,

\smallskip

 $\bullet\ J^{1,p}(\Omega):={\text{completion of }}{\mathscr{C}}_0(\Omega) {\text{ in the }} W^{1,p}(\Omega){\text{ norm}}$.

\medskip

For the properties of the  function spaces listed above and other spaces used in fluid dynamics, we refer to the monograph by G.P. Galdi \cite{Galdibook}.

We denote by $P$ the projector from $L^p$ onto $J^p$ and we call Stokes operator the composition $P\Delta$.

Throughout the paper, the symbols $c$, $c_1$, and $c_2$ denote positive constants whose numerical values are irrelevant for our purposes and may vary from line to line.

\bigskip
 
Let $g$ be a real valued function and consider the convolution
    \begin{equation}\label{convolution}
        f(t)=\int_a^{b}\frac{g(\tau)}{|t-\tau|^{1-\alpha}}{\mathrm d}\tau,\quad 0<\alpha\neq 1.
    \end{equation}
The following result is classical (see, {\it e.g.}, \cite{Oki,Stein}).
\begin{Lemma}\label{lemmaconvolution}
If $\displaystyle1<r<1/\alpha$, then \eqref{convolution} defines a continuous linear operator from $L^r(a,b)$ to $L^q(a,b)$, with $\frac
{1}{q}=\frac{1}{r}-\alpha$. In particular, there exists $c(r,q,\alpha)$ such that
\begin{equation}\label{stimaconvolution}
\|f\|_q\le c\|g\|_r.
\end{equation}
\end{Lemma}
We will also  use the following classical result (see \cite{CrispoMaremonti}).

\begin{Lemma}\label{lemmainterpolazione}
{\rm {(Gagliardo--Niremberg inequality in exterior domains).}}
Let $u\in L^q(\Omega)$ such that $D^m u\in L^r(\Omega)$, $1\le q< +\infty$ and $1\le r\le +\infty$. Then
\begin{equation}
    \|D^ju\|_p\le c\|D^m u\|_r^\theta\|u\|_q^{1-\theta},
\end{equation}
where $\frac{1}{p}=\frac{j}{2}+\theta\left(\frac{1}{r}-\frac{m}{2}\right)+(1-\theta)\frac{1}{q}$, for all $\theta\in[\frac{j}{m},1]$ $(\theta\ne 1$, if $1<r<+\infty$ and $m-j-2/r$ is a nonnegative integer$)$, and $c$ is a constant independent of $u$.
\end{Lemma}
The following result is proved in  \cite{Kato} for bounded domains and in \cite{Mare3} for exterior domains.
\begin{Lemma}\label{stimaconlaplaciano}
Let $u\in L^q(\Omega)$ such that $\Delta u\in J^p(\Omega)$, $\nabla u\in L_{\text{loc}}^p(\overline\Omega)$ and $\gamma_{\text{tr}}(u)=0\ on\ {\partial\Omega}$. Then
\begin{equation}
    \|u\|_\infty\le c\|P\Delta u\|_p^a\|u\|_q^{1-a},
\end{equation}
where $a\left(\frac{1}{p}-1\right)+(1-a)\frac{1}{q}=0$ and $c$ is a constant independent of $u$.
\end{Lemma}

\medskip 

We now recall some properties of the solutions to the linear system \eqref{Stokes}. First, we note that the result obtained in \cite{Farwig}, Theorem 2.4, for bounded domains,  still holds true for the other cases considered,  since it is based on the general interpolation result (\cite{Farwig}, Theorem 2.3) and on the following lemma. 

\begin{Lemma}\label{prepara}
    Let $1<r\le p\le\infty$ and $u_0\in  J^r(\Omega)$, $r\le\infty$. Then, the solution $(w,q)$ to \eqref{Stokes}--\eqref{condition8} satisfies
    \begin{equation}
        \|w(t)\|_p\le c t^{-\big(\frac{1}{r}-\frac{1}{p}\big)}\|u_0\|_r,\quad {\text{for }}t>0,
    \end{equation}
    where $c$ a constant independent of $u_0$.
\end{Lemma}
\begin{Proof}
    The result is stated in \cite{Farwig} for $\Omega$ bounded domain (see the references therein) and extends to the other cases under consideration thanks to the results in \cite{DanShi1, DanShi2, Mare5} for $r<\infty$ and \cite{Abe,Mare4} for $r=\infty$.\footnote{We recall that in $L^\infty(\Omega)$ the Helmholtz decomposition doesn't work. Hence, $J^\infty(\Omega)$ cannot be considered an element of the Helmholtz decomposition. Moreover, $J^\infty(\Omega)\subset{\mathscr {A}}$, where ${\mathscr{A}} =\{\phi\in L^\infty(\Omega):\ \int \phi\cdot\nabla h=0,\ \forall h\in L^1_{\text{loc}}(\Omega),\ \nabla h\in L^1(\Omega)\}$. The result in \cite{Mare4} holds for 
    $u_0\in {\mathscr{A}}$, as the ones proved in \cite{AbeGiga1,AbeGigaHieber}, that {\it a priori} hold locally in time.
}
\end{Proof}

Thanks to Lemma \ref{prepara} the following result holds.
 
 \begin{Lemma}\label{lemmaFarwigGiga}
     Let $(w,q)$ be the solution to \eqref{Stokes}--\eqref{condition8} with $u_0\in J^2(\Omega)$. Then
     \begin{equation}
    \displaystyle\int_0^{+\infty}\|w(\tau)\|_\infty^2{\text{d}}\tau \le c\|u_0\|_2^2,
\end{equation}
where $c$ is a constant independent of $u_0$.
 \end{Lemma}
 \begin{Proof}
 See Theorem 2.4. in \cite{Farwig}, p. 295.
 \end{Proof}

\bigskip 

Let $v=u-w$ and $\pi=p-q$, with $(u,p)$ and $(w,q)$ the solutions to \eqref{NavierStokes}--\eqref{condition2} and \eqref{Stokes}--\eqref{condition8}, respectively. The pair $(v,\pi)$ satisfies \eqref{forceterm}--\eqref{nullcondition2}
and takes into account the basic properties of $(u,p)$ and $(w,q)$, which are formulated
 in the following lemma.

\begin{Lemma}\label{proprieta}
    Let $u_0\in J^2(\Omega)$. Then, the solutions $(u,p)$ and $(w,q)$ to \eqref{NavierStokes}--\eqref{condition2} and \eqref{Stokes}--\eqref{condition8}, respectively, satisfy
      \begin{align}
        {}& u,w\in C([0,T); L^2(\Omega)),\label{uno}\\
        {}& u,w\in L^2(0,T;J^{1,2}(\Omega)),\label{due}\\
        {}& u,w\in L^2(\eta,T;W^{2,2}(\Omega)),\label{tre}
    \end{align}
    for all $T>\eta>0$.
\end{Lemma}
\begin{Proof}
    A classical result (cfr. \cite{OAL,DanShi1,DanShi2}) ensures the existence of  a solution $(w,q)$ to \eqref{Stokes}--\eqref{condition8} enjoying \eqref{uno}--\eqref{tre}. Concerning $u$, as a weak solution, it satisfies \eqref{uno} and \eqref{due} 
(cfr. \cite{OAL}); estimate \eqref{tre} holds as a consequence of its regularity for $t>0$.
 \end{Proof}

\medskip

For our purposes we also need to recall  
the special case of initial data in $L^1(\Omega)$ considered in \cite{Mare2}. Specifically, for $\varphi_0\in L^1(\Omega)$, 
the following initial--boundary value problem is analysed:
\begin{equation}\label{problemadatiL1}
\begin{array}{r@{}l}
\varphi_t-\Delta \varphi + \nabla Q =0&\qquad {\mathrm {in }}\ (0,+\infty)\times\Omega,\\
    \nabla \cdot \varphi= 0&\qquad {\mathrm {in }}\ (0,+\infty)\times\Omega,\\
    (\varphi(0),\psi)=(\varphi_0,\psi) &\quad\quad{\hbox {for any}}\ \psi\in\mathscr{C}_0(\Omega),\\
    \varphi=0 &\quad\quad{\hbox{on}}\ (0,+\infty)\times\partial\Omega.
\end{array}
\end{equation}
We remark that the initial condition \eqref{problemadatiL1}$_3$ means that  $\lim_{t\to 0} (\varphi(t)-\varphi_0,\psi)=0$, for any $\psi\in\mathscr{C}_0(\Omega)$. 
In \cite{Mare2} the following result is showed.

\begin{Lemma}\label{stimaconpeso}
Let $\varphi_0 \in L^1(\Omega)$. Then, there exists a unique solution $(\varphi,Q)$ to  \eqref{problemadatiL1} such that $\varphi\in C(\eta, T;J^q(\Omega))\cap L^\infty(\eta,T;J^{1,q}(\Omega)),\ D^2 \varphi,\nabla Q\in L^\infty(\eta,T;L^q(\Omega))$, for all $T>\eta>0$, $q>1$. Moreover,
\begin{equation}\label{semigroup1}
\|\varphi(t)\|_q\le \frac{c}{t^\mu}\|\varphi_0\|_1, \quad\forall t>0,
\end{equation}
where $\mu=1-\frac{1}{q}$,
    and
    \begin{equation}\label{semigroup2}
\|\nabla \varphi(t)\|_q\le \frac{c}{t^{\mu_1}}\|\varphi_0\|_1,\quad\mu_1=
\begin{cases}
    \frac{1}{2}+\mu\quad &{\hbox{if}}\ t\in(0,1],\\
    \frac{1}{2}+\mu\quad &{\hbox{if}}\ t>0\ {\hbox{and}}\ q\in(1,2],\\
    1\quad &{\hbox{if}}\ t>1\ {\hbox{and}}\ q>2.
\end{cases}
\end{equation}
In \eqref{semigroup1} and \eqref{semigroup2}, $q\in(1,\infty]$, $c$ is a constant independent of $\varphi_0$ and $\mu_1$ is sharp.

\smallskip

\noindent In particular, if $\varphi_0\in C^\infty_0(\Omega)$, then $\varphi\in \cap_{q>1}C([0,T);J^q(\Omega))$.

\smallskip

\end{Lemma}

\section{Proof of Theorem \ref{maintheorem}}\label{sectionproof}

\bigskip

{\sc Proof of Theorem \ref{maintheorem}} --
Let $u_0\in J^2(\Omega)$ and $(u,p)$ be the solution to \eqref{NavierStokes}--\eqref{condition2}. We represent $(u,p)$ by $(w+v,q+\pi)$, where $(w,q)$ and $(v,\pi)$ solve, respectively,
\begin{equation}\label{problemalineare}
\begin{array}{r@{}l}
    w_t-\Delta w + \nabla q &{}=0\qquad {\mathrm {in }}\ (0,+\infty)\times\Omega,\\
    \nabla \cdot w&{}= 0\qquad {\mathrm {in }}\ (0,+\infty)\times\Omega,\\
     w_{|t=0}&{}=u_0 \quad\;\,{\hbox {in}}\ \Omega,\\
    w&{}=0 \quad\quad{\hbox{on}}\ (0,+\infty)\times\partial\Omega,
\end{array}
\end{equation}
and 
\begin{equation}\label{problemanonlineare}
\begin{array}{r@{}l}
    v_t-\Delta v + \nabla \pi &{} =-u\cdot\nabla u\qquad{\mathrm {in }}\ (0,+\infty)\times\Omega,\\
    \nabla \cdot v &{} = 0\qquad\qquad\;\;\;\,{\mathrm {in }}\ (0,+\infty)\times\Omega,\\
    v_{|t=0}&{}=0 \qquad\qquad\;\;\;\,{\hbox {in}}\ \Omega,\\
    v&{}=0 \qquad\qquad\quad{\hbox{on}}\ (0,+\infty)\times\partial\Omega.%
\end{array}
\end{equation}

 By Lemma \ref{lemmaFarwigGiga}, 
$w$ satisfies
\begin{equation}\label{disu1}
    \displaystyle\int_0^{+\infty}\|w(t)\|_\infty^2 dt \le c\|u_0\|_2^2.
\end{equation}
As for $v$, 
by making use of Lemma \ref{proprieta} and Lemma \ref{stimaconpeso} we now establish that, for $p\in(4,+\infty)$,
\begin{equation}\label{stimaproposizione}
    \|v(t)\|_\infty\le c\int_0^t \frac{\|u(\tau)\|_p^2}{(t-\tau)^{\frac{1}{2}+\frac{2}{p}}}d\tau,\qquad \forall t>0,
    \end{equation}
    for some constant $c$.
   
 Indeed, let $\varphi_0\in C_0^\infty(\Omega)$ and $\varphi$ denote  the solution to problem \eqref{problemadatiL1}. Recalling that $\varphi(t-\tau,x)$ solves the backward problem in $(0,t)\times\Omega$,
  we multiply  \eqref{problemanonlineare}$_1$ by $\varphi(t-\tau,x)$ and integrate by parts, obtaining
    \begin{equation}\label{servecontinuita}
        (v(t),\varphi_0)=\int_0^t\int_\Omega u(\tau,x)\cdot\nabla\varphi(t-\tau,x)\cdot u(\tau,x)\ dx d\tau.
        \end{equation}
 In writing \eqref{servecontinuita}, we employed the continuity of $u$ and $\varphi$ in $L^2(\Omega)$, as guaranteed by Lemma \ref{proprieta} and Lemma  \ref{stimaconpeso}.
 Applying Hölder’s inequality yields
  \begin{equation}
      |(v(t),\varphi_0)|\le\int_0^t\|u(\tau)\|_p^2\|\nabla\varphi(t-\tau)\|_\frac{p}{p-2}d\tau.
  \end{equation}
Then, Lemma \ref{stimaconpeso} provides
\begin{equation}
      |(v(t),\varphi_0)|\le c\|\varphi_0\|_1\int_0^t\frac{\|u(\tau)\|_p^2}{(t-\tau)^{\frac{1}{2}+\frac{2}{p}}}d\tau,
  \end{equation}
  and \eqref{stimaproposizione} follows from the arbitrariness of $\varphi_0$ and
  the well--known formula 
  \begin{equation}
      \|v(t)\|_\infty= \text{sup}_{\varphi_0\in L^1(\Omega)}\frac{|(v(t),\varphi_0)|}{\|\varphi_0\|_1}.
  \end{equation}

With the estimate \eqref{stimaproposizione} just showed 
at hand,  
let us  observe that
\begin{equation}
    \int_0^t\frac{\|u(\tau)\|_p^2}{(t-\tau)^{\frac{1}{2}+\frac{2}{p}}}d\tau=\int_0^t\frac{\|u(\tau)\|_p^2}{(t-\tau)^{1-(\frac{1}{2}-\frac{2}{p})}}d\tau.
\end{equation}
Then, in virtue of  \eqref{stimaproposizione} and Lemma \ref{lemmaconvolution} with $\alpha=\frac{p-4}{2p}$, $q=2$ and $r=\frac{p}{p-2}$, we deduce
\begin{equation}\label{stimaintermedia}
\left(\displaystyle\int_0^{+\infty}\|v(t)\|_{\infty}^2 dt\right)^{1/2}
\le c \displaystyle\left(\int_0^{+\infty}\|u(t)\|_p^{2\cdot\frac{p}{p-2}}dt\right)^\frac{p-2}{p}.
\end{equation}
Applying the interpolation inequality (Lemma \ref{lemmainterpolazione}) with $j=0$, $m=1$, $q=r=2$ and $\theta=\frac{p-2}{p}$, we obtain
\begin{equation}
    \|u(t)\|_p\le c \|\nabla u(t)\|_2^{\frac{p-2}{p}}\|u(t)\|_2^{\frac{2}{p}}.
\end{equation}
Therefore, exploiting the energy identity \eqref{identity}, \eqref{stimaintermedia} implies
\begin{equation}
\begin{array}{{r@{}l}}
\Big(&\displaystyle\int_0^{+\infty}\|v(t)\|_{\infty}^2 dt\Big)^{1/2}\le c\left(\displaystyle\int_0^{+\infty}\|\nabla u(t)\|_2^2\|u(t)\|_2^{\frac{4}{p-2}}dt\right)^{\frac{p-2}{p}}\\[12pt]
{}&\le c\left(\displaystyle\int_0^{+\infty}\|\nabla u(t)\|_2^2\|u_0\|_2^{\frac{4}{p-2}}dt\right)^{\frac{p-2}{p}}=c\left(\|u_0\|_2^{\frac{4}{p-2}}\displaystyle\int_0^{+\infty}\|\nabla u(t)\|_2^2 dt\right)^{\frac{p-2}{p}}\\[12pt]
{}&\le c\left(\|u_0\|_2^{\frac{4}{p-2}}\|u_0\|_2^2\right)^{\frac{p-2}{p}}=c\|u_0\|_2^2.
\end{array}
\end{equation}
Hence, 
\begin{equation}\label{disu2}
    \displaystyle\int_0^{+\infty}\|v(t)\|_\infty^2 dt \le c\|u_0\|_2^4.
\end{equation}
Combining \eqref{disu1} with \eqref{disu2} and using the elementary inequality $|a+b|^2\le 2(|a|^2+|b|^2)$, for all $a,b\in{\mathbb R}$, we conclude that
    \begin{equation}
    \begin{array}{{r@{}l}}
    \displaystyle\int_0^{+\infty}\|u(t)\|_\infty^2 dt &\le 2 \displaystyle\int_0^{+\infty}\|w(t)\|_\infty^2 dt +2\int_0^{+\infty}\|v(t)\|_\infty^2 dt\\[12pt]
   {}& \le c_1\|u_0\|_2^2+ c_2\|u_0\|_2^4\le c\|u_0\|_2^2(1+\|u_0\|_2^2),
    \end{array}
\end{equation}
which completes the proof.\qed

\section{Proof of Theorem \ref{derivateseconde}}\label{sectionproof2}
{\sc Proof of Theorem \ref{derivateseconde}} --
Let $u_0\in J^2(\Omega)$ and $(u,p)$ the solution to \eqref{NavierStokes}--\eqref{condition2}. An integration by parts yields, $\forall t>0$,
\begin{equation}\label{parti}
\begin{array}{r@{}l}
     \displaystyle\frac{d}{dt}\displaystyle\int_\Omega|\nabla u|^2 dx &\ = 2\displaystyle\int_\Omega\nabla u_t\cdot\nabla u\:dx=-2\int_\Omega u_t\cdot\Delta u \:dx \\
     {}&=-2\displaystyle\int_\Omega|P\Delta u|^2dx+2\int_\Omega u\cdot\nabla u\cdot P\Delta u\:dx. 
\end{array}
\end{equation}
In virtue of H\"older's inequality and Young's inequality, the following  application of Lemma \ref{stimaconlaplaciano}
\begin{equation}\label{stimalapl}
    \|u\|_\infty\le c\|P\Delta u\|_2^{1/2}\|u\|_2^{1/2},\qquad\forall t>0,
\end{equation}
provides
\begin{equation}
    \int_\Omega u\cdot\nabla u\cdot P\Delta u\:dx\le c\|u\|_2^{1/2}\|\nabla u\|_2\|P\Delta u\|_2^{3/2}\le c\|u\|_2^2\|\nabla u\|_2^4 + \varepsilon \|P\Delta u\|_2^2,
\end{equation}
with $\varepsilon$ arbitrarily small.
Thus, \eqref{parti} and \eqref{identity} yield
\begin{equation}\label{stimaderivatagradiente}
     \frac{d}{dt}\|\nabla u(t)\|_2^2+\|P\Delta u(t)\|_2^2\le c\|u_0\|_2^2\|\nabla u(t)\|_2^4,\quad \forall t>0.
\end{equation}
By virtue of \eqref{identity} once again, 
\begin{equation}\label{stimauniformegradiente}
\int_0^t\|\nabla u(\tau)\|_2^2 d\tau\le \|u_0\|_2^2.
\end{equation}
Moreover,
\begin{equation}
\begin{array}{r@{}l}
    \displaystyle\int_0^t\|\nabla u(\tau)\|_2^2 d\tau &\ =\displaystyle\int_0^t\frac{d}{d\tau}(\tau\|\nabla  u(\tau)\|_2^2)d\tau-\int_0^t\tau\frac{d}{d\tau}\|\nabla u(\tau)\|_2^2d\tau\\
    &=t\|\nabla u(t)\|_2^2-\displaystyle\int_0^t\tau\frac{d}{d\tau}\|\nabla u(\tau)\|_2^2 d\tau.
\end{array}    
\end{equation}
Therefore, from \eqref{stimaderivatagradiente} and \eqref{stimauniformegradiente} we obtain
\begin{equation}\label{stimaperzeta}
    t\|\nabla u(t)\|_2^2\le \|u_0\|_2^2+c\|u_0\|_2^2\int_0^t\tau\|\nabla u(\tau)\|_2^4d\tau.
\end{equation}
Setting $z(t):=t\|\nabla u(t)\|_2^2$, \eqref{stimaperzeta} writes
\begin{equation}\label{gron}
    z(t)\le \|u_0\|_2^2+c\|u_0\|_2^2\int_0^tz(\tau)\|\nabla u(\tau)\|_2^2d\tau.
\end{equation}
In virtue of the Gronwall lemma and \eqref{stimauniformegradiente}, \eqref{gron} implies
\begin{equation}
    z(t)\le\|u_0\|_2^2 e^{c\|u_0\|_2^2\int_0^t\|\nabla u(\tau)\|_2^2d\tau}\le\|u_0\|_2^2 e^{c\|u_0\|_2^4},
\end{equation}
that is,
\begin{equation}\label{peso}
\|\nabla u (t)\|_2\le\frac{\|u_0\|_2 e^{c\|u_0\|_2^4}}{t^{1/2}}.
\end{equation}
Now, integrating \eqref{stimaderivatagradiente} on $(s,t)$, $s>0$, we get
\begin{equation}\label{integrato}
    \|\nabla u(t)\|_2^2+\displaystyle\int_s^t\|P\Delta u(\tau)\|_2^2d\tau
    \le\|\nabla u(s)\|_2^2+c\|u_0\|_2^2\displaystyle\int_s^t\|\nabla u(\tau)\|_2^4d\tau.
\end{equation}
Then, integrating
\eqref{integrato} with respect to $s$ on $(0,t)$, we obtain
\begin{equation}\label{pippo}
    \int_0^tds\int_s^t\|P\Delta u(\tau)\|_2^2d\tau\le\int_0^t\|\nabla u(s)\|_2^2ds+c\|u_0\|_2^2\int_0^tds\int_s^t\|\nabla u(\tau)\|_2^4d\tau.
\end{equation}
Taking into account \eqref{comportamentoneltempo}, an integration by parts yields
\begin{equation}
\begin{array}{r@{}l}
\displaystyle\int_0^t ds\int_s^t \|P\Delta u(\tau)\|_2^2 &d\tau =-\displaystyle\int_0^t s\frac{d}{ds}\int_s^t \|P\Delta u(\tau)\|_2^2 d\tau\\
{}&=\displaystyle\int_0^t s\|P\Delta u(s)\|_2^2ds.
\end{array}
\end{equation}
Analogously,
\begin{equation}
\displaystyle\int_0^t ds\int_s^t \|\nabla u(\tau)\|_2^4 d\tau =\displaystyle\int_0^t s\|\nabla u(s)\|_2^4ds.
\end{equation}
Therefore, recalling \eqref{stimauniformegradiente}, \eqref{pippo} implies
\begin{equation}\label{bo}
    \int_0^t s\|P\Delta u(s)\|_2^2ds\le\|u_0\|_2^2+c\|u_0\|_2^2\int_0^t s\|\nabla u(s)\|_2^4ds.
\end{equation}
By applying \eqref{peso} and \eqref{stimauniformegradiente}, from \eqref{bo} we get
\begin{equation}\label{ultima}
   \int_0^t s\|P\Delta u(s)\|_2^2ds\le c \|u_0\|_2^2(1+\|u_0\|_2^4 e^{c_1\|u_0\|_2^4}).
\end{equation}
Estimate \eqref{ultima} holds with $c$ and $c_1$ constants independent of $t>0$, so that, letting $t\to +\infty $ yields \eqref{primastimateorema2}.
Now, in virtue of \eqref{stimalapl} and \eqref{identity}, one has
\begin{equation}\label{stimalapl2}
    \|u\|_\infty\le c\|P\Delta u\|_2^{1/2}\|u_0\|_2^{1/2},\qquad\forall t>0.
\end{equation}
Hence, \eqref{primastimateorema2} implies \eqref{secondastimateorema2}.
\qed

\bigskip
\bigskip

{\bf Acknowledgements.} This research was performed under the auspices
of GNFM Gruppo Nazionale per la Fisica Matematica -- INdAM Istituto
Nazio\-na\-le di Alta Matematica ``Francesco Severi''.

\end{document}